    \documentclass[12pt]{article}
    \usepackage{t1enc}
    \usepackage[latin1]{inputenc}
    \usepackage[english]{babel}
        \usepackage{latexsym}
        \newcommand{\bp}{{\bf P}}

\newcommand{\be}{{\bf E}}

\def\eqe{\stackrel{{\cal D}}{=}}
\newtheorem{theorem}{Theorem}[section]

\newtheorem{lemma}{Lemma}[section]

\def\begg{begin{equation}}
\def\endd{\end{equation}}

\def\z2{{\cal Z}_2}

\def\begg{\begin{equation}}
\def\endd{\end{equation}}
\def\bege{\begin{eqnarray}}
\def\ende{\end{eqnarray}}

\def\pe{{\bf P}}

\begin{document}

\centerline{\Large\bf TRANSIENT NEAREST NEIGHBOR RANDOM WALK}
\medskip
\centerline{\Large\bf AND BESSEL PROCESS}
\medskip

\bigskip \bigskip \bigskip \bigskip \bigskip

\renewcommand{\thefootnote}{1} \noindent
\textbf{Endre Cs\'{a}ki}\footnote{Research supported by the
Hungarian National Foundation for Scientif\/ic Research, Grant No.
  K 61052 and K 67961.}\newline Alfr\'ed R\'enyi Institute of
Mathematics, Hungarian Academy of Sciences, Budapest, P.O.B. 127,
H-1364, Hungary. E-mail address: csaki@renyi.hu

\bigskip

\renewcommand{\thefootnote}{2} \noindent \textbf{Ant\'{o}nia
F\"{o}ldes}\footnote{Research supported by a PSC CUNY Grant, No.
69020-0038.}\newline Department of Mathematics, College of Staten
Island, CUNY, 2800 Victory Blvd., Staten Island, New York 10314,
U.S.A. E-mail address: foldes@mail.csi.cuny.edu

\bigskip

\noindent \textbf{P\'al R\'ev\'esz}$^1$ \newline Institut f\"ur
Statistik und Wahrscheinlichkeitstheorie, Technische Universit\"at
Wien, Wiedner Hauptstrasse 8-10/107 A-1040 Vienna, Austria. E-mail
address: reveszp@renyi.hu

\bigskip \bigskip \bigskip

\noindent \textit{Abstract:} We prove strong invariance principle
between a transient Bessel process and a certain nearest neighbor
(NN) random walk that is constructed from the former by using
stopping times. It is also  shown that their local times are close
enough to share the same strong limit theorems. It is shown
furthermore, that if the difference between the distributions of two
NN random walks are small, then the walks themselves can be
constructed so that they are close enough. Finally, some
consequences concerning strong limit theorems are discussed.

\bigskip


\bigskip

 \vspace{.1cm}

\noindent Running head: NN random walk and Bessel process.

\vfill
\renewcommand{\thesection}{\arabic{section}.}
\section{Introduction}

\renewcommand{\thesection}{\arabic{section}} \setcounter{equation}{0}
\setcounter{theorem}{0} \setcounter{lemma}{0}

In this paper we consider a nearest neighbor (NN) random walk,
defined as follows: let $X_0=0,\ X_1,X_2,\ldots$ be a Markov
chain with
\begin{eqnarray}\label{defff}
E_i&:=&\bp(X_{n+1}=i+1\mid X_n=i)=1-\bp(X_{n+1}=i-1\mid X_n=i)\\
&=&\left\{\begin{array}{ll} 1\quad & {\rm if}\quad  i=0\\
\nonumber
 1/2+p_i\quad & {\rm if}\quad i=1,2,\ldots,
\end{array}\right.
\end{eqnarray}
where $-1/2\leq p_i\leq 1/2,\ i=1,2,\ldots$. In case $0< p_i\leq
1/2$ the sequence $\{X_i\}$ describes the motion of a particle which
starts at zero, moves over the nonnegative integers and going away
from 0 with a larger probability than to the direction of 0. We will
be interested in the case when $p_i\sim B/4i$ with $B>0$ as
$i\to\infty$. We want to show that in certain sense, this Markov
chain is a discrete analogue of continuous Bessel process and
establish a strong invariance principle between these two processes.

The properties of the discrete model, often called birth and death
chain, connections with orthogonal polynomials in particular, has
been treated extensively in the literature. See e.g. the
classical paper by Karlin and McGregor \cite{KMG}, or more recent
papers by Coolen-Schrijner and Van Doorn \cite{C-SD} and Dette
\cite{DE01}. In an earlier paper \cite{CSFR} we investigated the
local time of this Markov chain in the transient case.

There is a well-known  result in the literature (cf. e.g. Chung
\cite{CH}) characterizing those sequences $\{p_i\}$ for which
$\{X_i\}$ is transient (resp. recurrent).

\smallskip

\noindent
 {\bf Theorem A:} (\cite{CH}, page 74) {\it Let $X_n$ be a Markov chain
 with transition probabilities given in {\rm (\ref{defff})} with
 $-1/2<p_i<1/2$, $i=1,2,\ldots$ Define

\begg U_i:={\frac{1-E_i}{E_i}}={\frac{1/2-p_i}{1/2+p_i}}
\label{uif}
\endd
Then $X_n$  is transient if and only if}
$$\sum_{k=1}^{\infty} \prod_{i=1}^k U_i < \infty.$$

As a consequence, the Markov chain $(X_n)$ with $p_R\sim B/4R,\,
R\to\infty$ is transient if $B>1$ and recurrent if $B<1$.

The Bessel process of order $\nu$, denoted by $Y_\nu(t),\, t\geq
0$ is a diffusion process on the line with generator
$$
\frac{1}{2}\frac{d^2}{dx^2}+\frac{2\nu+1}{2x}\frac{d}{dx}.
$$
$d=2\nu+2$ is the dimension of the Bessel process. If $d$ is a
positive integer, then $Y_\nu(\cdot)$ is the absolute value of a
$d$-dimensional Brownian motion. The Bessel process $Y_\nu(t)$ is
transient if and only if $\nu>0$.

The properties of the Bessel process were extensively studied in
the literature. Cf. Borodin and Salminen \cite{BS}, Revuz and Yor
\cite{RY}, Knight \cite{KN81}.

Lamperti \cite{LA63} determined the limiting distribution of $X_n$ and
also proved a weak convergence theorem in a more general setting. His
result in our case reads as follows.

\smallskip
\noindent {\bf Theorem B:} (\cite{LA63}) {\it Let $X_n$ be a Markov
chain with transition probabilities given in {\rm (\ref{defff})}
with  $-1/2<p_i<1/2$, $i=1,2,\ldots$ If
$\,\,\lim_{R\to\infty}\,Rp_R=B/4>-1/4$,} {\it then the following
weak convergence holds}:
$$
\frac{X_{[nt]}}{\sqrt{n}}\Longrightarrow Y_{(B-1)/2}(t)
$$
{ \it in the space  D}[0,1]. {\it In particular,}
$$
\lim_{n\to \infty} \bp\left(\frac{X_n}{\sqrt{n}}<x\right)=
\frac{1}{2^{B/2-1/2}\Gamma(B/2+1/2)}\int_0^x u^Be^{-u^2/2}\, du.
$$

 In Theorems A and B values of $p_i$ can be negative.
In the sequel however we deal only with the case when $p_i$ are
non-negative, and the chain is transient, which will be assumed
throughout without mentioning it.

Let
\begg
D(R,\infty):=1+\sum_{j=1}^{\infty}\prod_{i=1}^j U_{R+i},
\endd
and define
\begg
p_R^*:=\frac{\frac{1}{2}+p_R}{D(R,
\infty)}=1-q^*_R
\endd
Now let $\xi(R,\infty)$, $R=0,1,2,\ldots$ be the total local time
at $R$ of the Markov chain $\{X_n\}$, i.e. \begg
\xi(R,\infty):=\#\{n\geq 0: X_n=R\}. \label{loct}
\endd
\noindent {\bf Theorem C:} (\cite{CSFR}) {\it For a transient NN
random walk}
\begg \pe(\xi(R,\infty)=k)=p_R^*(q_R ^*)^{k-1},
\qquad k=1,2,\ldots \label{locelo}
\endd
Moreover, $\eta(R,t)$, $R>0$ will denote the local time of the
Bessel process, i.e.
$$
\eta(R,t):=\lim_{\varepsilon\to 0}\frac1{2\varepsilon}\int_0^t
I\{Y_\nu(s)\in (R-\varepsilon,R+\varepsilon)\}\, ds,\qquad
\eta(R,\infty):=\lim_{t \to \infty }\eta(R,t).
$$
It is well-known that $\eta(R,\infty)$ has exponential
distribution (see e.g. \cite{BS}).
 \begg \pe(\eta(R,\infty)< x)=
1-\exp\left(-\frac{\nu }{R}\,x\right).
\endd
\newline
For $0<a<b$ let
\begin{equation}
\tau:=\tau(a,b)=\min\{t\geq 0:\, Y_\nu(t)\notin (a,b)\}.
\label{stopping}
\end{equation}
Then we have (cf. Borodin and Salminen \cite{BS}, Section 6,
3.0.1 and 3.0.4).

\smallskip\noindent {\bf Theorem D:} {\it For} $0<a<x<b$ {\it we have}
\begin{equation}
\bp_x(Y_\nu(\tau)=a)=1-\bp_x(Y_\nu(\tau)=b)=
\frac{x^{-2\nu}-b^{-2\nu}}{a^{-2\nu}-b^{-2\nu}}, \label{hit}
\end{equation}
\begin{equation}
\be_xe^{-\alpha\tau}=\frac{S_\nu(b\sqrt{2\alpha},x\sqrt{2\alpha})+
S_\nu(x\sqrt{2\alpha},a\sqrt{2\alpha})}
{S_\nu(b\sqrt{2\alpha},a\sqrt{2\alpha})}, \label{lap}
\end{equation}

\noindent {\it where}
\begin{equation}
S_\nu(u,v)=(uv)^{-\nu}(I_\nu(u)K_\nu(v)-K_\nu(u)I_\nu(v)),
\label{snu}
\end{equation}

\noindent $I_\nu$ {\it and} $K_\nu$ {\it being the modified Bessel
functions of the first and second kind, resp.}

Here and in what follows $\bp_x$
and $\be_x$ denote conditional probability, resp. expectation
under $Y_{\nu}(0)=x.$ For simplicity we  will use $\bp_0=\bp,$ and
$\be_0=\be.$

 Now consider
$Y_\nu(t),\, t\geq 0$, a Bessel process of order $\nu$,
$Y_\nu(0)=0$, and let $X_n,\, n=0,1,2,\ldots$ be an NN random
walk with  $p_0=p_1=1/2$,
\begin{equation}
p_R=\frac{(R-1)^{-2\nu}-R^{-2\nu}}{(R-1)^{-2\nu}-(R+1)^{-2\nu}}-\frac12,\qquad
R=2,3,\ldots \label{pr}
\end{equation}

Our main results are strong invariance principles concerning Bessel
process, NN random walk and their local times.

\begin{theorem}
On a suitable probability space we can construct a Bessel process
$\{Y_\nu(t),\, t\geq 0\},$ $\nu>0$ and an NN random walk $\{X_n,\,
n=0,1,2,\ldots\}$ with $p_R$ as in {\rm (\ref{pr})} such that for
any $\varepsilon>0$, as $n\to\infty$ we have
\begin{equation}
Y_\nu(n)-X_n=O(n^{1/4+\varepsilon})\qquad {\rm a.s.} \label{inv}
\end{equation}
\end{theorem}

Our strong invariance principle for local times reads as follows.

\begin{theorem} Let $Y_\nu(t)$ and $X_n$ as in {\rm Theorem 1.1}
and let $\eta$ and $\xi$ their respective local times. As
$R\to\infty$, we have
\begin{equation}
\xi(R,\infty)-\eta(R,\infty)=O(R^{1/2}\log R)\quad {\rm a.s.}
\end{equation}
\end{theorem}

We prove the following strong invariance principle  between two NN
random walks.

\begin{theorem}
Let $\{X_n^{(1)}\}_{n=0}^\infty$ and $\{X_n^{(2)}\}_{n=0}^\infty$
be two NN random walk with $p_j^{(1)}$ and $p_j^{(2)}$, resp.
Assume that
\begin{equation}
\left|p_j^{(1)}-\frac{B}{4j}\right|\leq \frac{C}{j^\gamma}
\label{pj1}
\end{equation}
and
\begin{equation}
\left|p_j^{(2)}-\frac{B}{4j}\right|\leq \frac{C}{j^\gamma}
\label{pj2}
\end{equation}
$j=1,2,\ldots$ with $B>1$, $1<\gamma\leq 2$ and some non-negative
constant $C$. Then on a suitable probability space one
can construct $\{X_n^{(1)}\}$ and $\{X_n^{(2)}\}$ such that as $n\to
\infty$
$$
|X_n^{(1)}-X_n^{(2)}|=O((X_n^{(1)}+X_n^{(2)})^{2-\gamma})
=O((n\log\log n)^{1-\gamma/2})\quad{\rm a.s.}
$$
\end{theorem}

\medskip

The organization of the paper is as follows. In Section 2 we will
present some well-known facts and prove some preliminary results.
Sections 3-5 contain the proofs of Theorems 1.1-1.3, respectively.
In Section 6 we prove strong theorems (most of them are integral
tests) which easily follow from Theorems 1.1 and 1.2 and the
corresponding results for Bessel process. In Section 7, using our
Theorem 1.3 in both directions, we prove an integral
test for the local time of the NN-walk, and a strong theorem for the
speed of escape of the Bessel process.

\renewcommand{\thesection}{\arabic{section}.}
\section{Preliminaries}
\renewcommand{\thesection}{\arabic{section}} \setcounter{equation}{0}
\setcounter{theorem}{0} \setcounter{lemma}{0}

\begin{lemma}
Let $Y_\nu(\cdot)$ be a Bessel process starting from $x=R$ and let
$\tau$ be the stopping time defined by {\rm (\ref{stopping})} with
 $a=R-1$ and $b=R+1$. Let $p_R$ be defined  by {\rm (\ref{pr})}.
Then as $R\to\infty$
\begin{equation}
p_R=\frac{2\nu+1}{4R}+O\left(\frac1{R^2}\right), \label{pr2}
\end{equation}
\begin{equation}
\be_R(\tau)=1+O\left(\frac{1}{R}\right), \label{ertau}
\end{equation}
\begin{equation}
Var_R(\tau)=O(1). \label{vartau}
\end{equation}
\end{lemma}

\vspace{2ex}\noindent {\bf Proof:} For $\nu=1/2$, i.e. for
$d=3$-dimensional Bessel process, in case $x=R$, $a=R-1$, $b=R+1$
we have
$$
\be_R(e^{\lambda\tau})=\frac1{\cos(\sqrt{2\lambda})}
$$
which does not depend on $R$. We prove that this holds
asymptotically in general, when $\nu>0$.

Using the identity (cf. \cite{BS}, page 449 and \cite{WA}, page 78)

\begin{eqnarray}\label{deef}
K_\nu(x)&=&
\left\{\begin{array}{ll}\displaystyle{\frac{\pi}{2\sin(\nu\pi)}
(I_{-\nu}(x)-I_\nu(x))}\quad
{\rm if} \,\, \nu \,\,{\rm is\,\, not\,\, an\,\, integer\,\,}\\  & \\
\nonumber
  \lim_{\mu\to \nu}K_{\mu}(x) \quad{\rm if\,\,} \nu\,\, {\rm is\,\,
  an\,\,
  integer}
\end{array}\right.
\end{eqnarray}

\noindent and the series expansion
$$
I_\nu(x)=\sum_{k=0}^\infty\frac{(x/2)^{\nu+2k}}{k!\Gamma(\nu+k+1)},
$$
one can see that the coefficient of $-\alpha$ in the  Taylor series
expansion of  the Laplace transform (\ref{lap}) is
$$
\be_x(\tau)=\frac1{2(\nu+1)}\frac{(b^2-x^2)a^{-2\nu}+(x^2-a^2)b^{-2\nu}
-(b^2-a^2)x^{-2\nu}}{a^{-2\nu}-b^{-2\nu}}
$$
from which by putting $x=R$, $a=R-1$, $b=R+1$, we obtain
$$
\be_R(\tau)=\frac1{2(\nu+1)}\frac{(2R+1)(R-1)^{-2\nu}+(2R-1)(R+1)^{-2\nu}
-4R^{1-2\nu}}{(R-1)^{-2\nu}-(R+1)^{-2\nu}}
$$
giving (\ref{ertau}) after some calculations.

(\ref{vartau}) can also be obtained similarly, but it seems quite
complicated. A simpler argument is to use moment generating
function and expansion of the Bessel functions for imaginary
arguments near infinity. Put $\alpha=-\lambda$ into (\ref{lap})
to obtain
\begin{equation}
\be_x(e^{\lambda\tau})=\frac{S_\nu(ib\sqrt{2\lambda},ix\sqrt{2\lambda})+
S_\nu(ix\sqrt{2\lambda},ia\sqrt{2\lambda})}
{S_\nu(ib\sqrt{2\lambda},ia\sqrt{2\lambda})}, \label{mgf}
\end{equation}
where $i=\sqrt{-1}$. We use the following  asymptotic expansions
(cf. Erd\'elyi et al.\,\cite{EA}, page 86, or Watson \cite{WA}, pages
202, 219)

$$
I_\nu(z)=(2\pi
z)^{-1/2}\left(e^z+ie^{-z+i\nu\pi}+O(|z|^{-1})\right),
$$
$$
K_\nu(z)=\left(\frac{\pi}{2z}\right)^{1/2}\left(e^{-z}+O(|z|^{-1})\right).
$$

 Hence one obtains for  $\lambda>0$ fixed, and $x<b,\,$
$$
S_\nu(ib\sqrt{2\lambda},ix\sqrt{2\lambda})= (-2\lambda b x)^{-\nu}
(I_\nu(ib\sqrt{2\lambda})K_{\nu}(ix\sqrt{2\lambda})
-I_\nu(ix\sqrt{2\lambda})K_{\nu}(ib\sqrt{2\lambda}))
$$
$$
=\frac12 (-2\lambda b x)^{-\nu-1/2}
\left(e^{i(b-x)\sqrt{2\lambda}}-e^{-i(b-x)\sqrt{2\lambda}}
+O\left(\frac1x\right)\right), \quad x\to\infty.
$$
One can obtain asymptotic expansions similarly for
$S_\nu(ix\sqrt{2\lambda},ia\sqrt{2\lambda})$,
$S_\nu(ib\sqrt{2\lambda},ia\sqrt{2\lambda})$. Putting these into
(\ref{mgf}), with $x=R$, $a=R-1$, $b=R+1$, we get as $R \to
\infty$
$$
\be_R(e^{\lambda\tau})=\frac{(R^2+R)^{-\nu-1/2}+(R^2-R)^{-\nu-1/2}}
{(R^2-1)^{-\nu-1/2}}\, \,
\frac{e^{i\sqrt{2\lambda}}-e^{-i\sqrt{2\lambda}}
+O\left(\frac1R\right)}
{e^{2i\sqrt{2\lambda}}-e^{-2i\sqrt{2\lambda}}
+O\left(\frac1R\right)}
$$
$$
=\frac1{\cos(\sqrt{2\lambda})}+O\left(\frac1R\right).
$$
Hence putting $\lambda=1$, there exists a constant $C$ such that
$\be_R(e^{\tau})\leq C$ for all $R=1,2,\ldots$ By Markov's
inequality we have
$$
\bp_R(\tau>t)=\bp_R(e^\tau>e^t)\leq Ce^{-t},
$$
from which $\be_R(\tau^2)\leq 2C$, implying (\ref{vartau}). $\Box$

Here and throughout $C,C_1,C_2,\ldots$ denotes unimportant positive
(possibly random) constants whose values may change from line to line.

Recall the definition of the upper and lower classes for a
stochastic process $Z(t),\, t\geq 0$ defined on a probability space
$(\Omega,{\cal F}, P)$  (cf. R\'ev\'esz \cite{R05}, p. 33).

\bigskip

The function $a_1(t)$ belongs to the {\it upper-upper class} of
$Z(t)$ ($a_1(t)\in {\rm UUC}(Z(t)$) if  for almost all $\omega \in
\Omega $ there exists a $t_0(\omega)>0 $ such that
 $Z(t)<a_1(t)$ if  $t>t_0(\omega).$

\bigskip

The function $a_2(t)$ belongs to the {\it upper-lower class} of
$Z(t)$ ($a_1(t)\in {\rm ULC}(Z(t)$) if for almost all $\omega \in
\Omega $  there exists a sequence of positive numbers
$0<t_1=t_1(\omega)<t_2=t_2(\omega)< \ldots $ with
$\lim_{i\to\infty}t_i=\infty$ such that $Z(t_i)\geq a_2(t_i)$,
$(i=1,2,\ldots).$

\bigskip

The function $a_3(t)$ belongs to the {\it lower-upper class } of
$Z(t)$ ($a_3(t)\in {\rm LUC}(Z(t)$) if for almost all $\omega \in
\Omega $  there exists a sequence of positive numbers
$0<t_1=t_1(\omega)<t_2=t_2(\omega)< \ldots $ with
$\lim_{i\to\infty}t_i=\infty$ such that $Z(t_i)\leq a_3(t_i)$,
$(i=1,2,\ldots).$

\bigskip

The function $a_4(t)$ belongs to the {\it lower-lower class} of
$Z(t)$ ($a_4(t)\in {\rm LLC}(Z(t)$) if for almost all $\omega \in
\Omega $ there exists a $t_0(\omega)>0 $ such that
 $Z(t)> a_4(t)$ if  $t>t_0(\omega).$ \bigskip

The following lower class results are due to Dvoretzky
and Erd\H os \cite{DE} for integer $d=2\nu+2$. In the general
case when $\nu>0$, the proof is similar (cf. also Knight
\cite{KN81} and Chaumont and Pardo \cite{CP} in the case of
positive self-similar Markov processes).

\bigskip
\noindent{\bf Theorem E:} {\it Let $\nu>0$ and let $b(t)$ be a
non-increasing, non-negative function.}

\begin{itemize}
\item $t^{1/2}b(t) \in {\rm LLC}(Y_\nu(t)) $ \qquad {\it if and only if}
\qquad
  $\displaystyle{
\int_1^\infty (b(2^t))^{2\nu}\, dt< \infty.}$ \label{at}
\end{itemize}
\medskip
It follows e.g. that in case $\nu>0$, for any $\varepsilon>0$ we
have
\begin{equation}
Y_\nu(t)\geq t^{1/2-\varepsilon} \label{ylower}
\end{equation}
almost surely for all sufficiently large $t$.

In fact, from our invariance principle it will follow that the
integral test in Theorem E holds also for our Markov chain
$(X_n)$. In the proof however we need an analogue of
(\ref{ylower}) for $X_n$.

 One can easily calculate the exact distribution of
$\xi(R,\infty),$ the total local time of $X_n$ of Theorem 1.1
according to Theorem C.
\newline
 {\bf Lemma A:} {\it If $p_R$ is given by {\rm (\ref{pr})},
then} $\xi(R,\infty)$ {\it has geometric distribution {\rm
(\ref{locelo})} with}
\begg
p_R^*=\frac{\frac{1}{2}+p_R}{D(R,\infty)}=
\frac{(\frac{1}{2}+p_R)((R+1)^{2\nu}-R^{2\nu})}{(R+1)^{2\nu}}=
\frac{\nu}{R}+O\left(\frac{1}{R}\right).
\endd

\begin{lemma} For any $\delta>0$ we have
$$
X_n\geq n^{1/2-\delta}
$$
almost surely for all large enough $n$.
\end{lemma}

\noindent{\bf Proof:} From Lemma A it is easy to
conclude that almost surely for some $R_0>0$
$$\xi(R,\infty)\leq CR\log R$$
if $R\geq R_0$, with some random positive constant $C$. Hence the
time $\displaystyle{\sum_{R=1}^S}\xi(R,\infty)$ which the particle
spent up to $\infty$ in $[1,S]$ is less than
$$\sum_{R=1}^{R_0-1}\xi(R,\infty)+C\sum_{R=R_0}^S R\log R\leq
C_1S^{2+\delta}$$ with some (random) $C_1>0$. Consequently, after
$C_1S^{2+\delta}$ steps the particle will be farther away from
the origin than $S$. Let
$$n=[C_1S^{2+\delta}],$$
then
$$S\geq\left(\frac{n}{C_1}\right)^{1/(2+\delta)}$$
and hence
$$
X_n\geq \left(\frac{n}{C_1}\right)^{1/(2+\delta)}\geq
n^{1/2-\delta}
$$
for $n$ large enough. This proves the Lemma. $\Box$


\renewcommand{\thesection}{\arabic{section}.}
\section{Proof of Theorem 1.1}
\renewcommand{\thesection}{\arabic{section}} \setcounter{equation}{0}
\setcounter{theorem}{0} \setcounter{lemma}{0}

Def\/ine the sequences $(\tau_n)$, $t_0=0,$
$t_n:=\tau_1+\ldots+\tau_n$ as follows:
\begin{eqnarray*}
\tau_1&:=&\min\{t:\ t>0,\ Y_\nu(t)=1\},\\
\tau_2&:=&\min\{t:\ t>0,\ Y_\nu(t+t_1)=2\},\\
\tau_n&:=&\min\{t:\ t>0,\
|Y_\nu(t+t_{n-1})-Y_\nu(t_{n-1})|=1\}\quad {\rm for} \quad
n=3,4,\ldots
\end{eqnarray*}

Let $X_n=Y_\nu(t_n)$. Then (cf. (\ref{pr})) it is an NN random
walk with $p_0=p_1=1/2$,
$$
p_R=\frac{(R-1)^{-2\nu}-R^{-2\nu}}{(R-1)^{-2\nu}-(R+1)^{-2\nu}}-\frac12,
\qquad R=2,3,\ldots
$$

Let ${\cal F}_n$ be the $\sigma$-algebra generated by $(\tau_k,\
Y_\nu(\tau_k))_{k=1}^n$ and consider
$$
M_n:=\sum_{i=1}^n(\tau_i-\be(\tau_i\mid {\cal F}_{i-1})).
$$
Then the sequence $(M_n)_{n\geq 1}$ is a martingale with respect to
$({\cal F}_n)_{n\geq 1}$. It follows from (\ref{ertau}) of Lemma 2.1
that for $i=2,3,\ldots$ we have
$$
\be(\tau_i\mid {\cal F}_{i-1})=\be(\tau_i\mid Y_\nu(t_{i-1}))
=1+O\left(\frac{1}{Y_\nu(t_{i-1})}\right).
$$
Hence
$$
|t_n-n|\leq |M_n|
+|\tau_1-1|+C_1\sum_{i=2}^n\frac{1}{Y_\nu(t_{i-1})}= |M_n|
+|\tau_1-1|+C_1\sum_{i=2}^n\frac{1}{X_{i-1}}
$$
with some (random) constant $C_1$. By (\ref{vartau}) of Lemma 2.1
we have $\be M_n^2\leq Cn$. Let $\varepsilon>0$ be arbitrary and
define $n_k=[k^{1/\varepsilon}]$. From the martingale inequality we
get
$$
\bp\left(\max_{n_{k-1}\leq n\leq n_k}|M_n|\geq
C_1n_{k-1}^{1/2+\varepsilon}\right) \leq
\frac{C_2}{n_k^{2\varepsilon}},
$$
hence we obtain by Borel-Cantelli lemma
$$
\max_{n_{k-1}\leq n\leq n_k}|M_n|\leq C_1
n_{k-1}^{1/2+\varepsilon}
$$
almost surely for large $k$. Hence we also have
$$
|M_n|=O(n^{1/2+\varepsilon})\qquad {\rm a.s.}
$$
By Lemma 2.2
$$
\sum_{i=2}^n\frac1{X_{i-1}}=O(n^{1/2+\varepsilon})\qquad {\rm
a.s.},
$$
consequently
\begin{equation}
|t_n-n|=O(n^{1/2+\varepsilon})\qquad{\rm a.s.} \label{tminusn}
\end{equation}
It is well-known (cf. \cite{BS}, p. 69) that $Y_\nu(t)$ satisfies
the stochastic differential equation \begg
dY_\nu(t)=dW(t)+\frac{2\nu+1}{2Y_\nu(t)}dt, \label{dife}
\endd
where $W(t)$ is a standard Wiener process. Hence
$$
X_n-Y_\nu(n)=Y_\nu(t_n)-Y_\nu(n)=
W(t_n)-W(n)+\int_{t_n}^n\frac{2\nu+1}{2Y_\nu(s)}\, ds,
$$
consequently,
$$
|X_n-Y_\nu(n)|\leq |W(t_n)-W(n)|+\frac{(2\nu+1)|t_n-n|}{2}
\max_{\min(n,t_n)\leq t\leq \max(n,t_n)}\frac1{Y_\nu(t)}.
$$

Now by (\ref{tminusn}) and (\ref{ylower}) the last term is
$O(n^{2\varepsilon})$ almost surely and since for the increments
of the Wiener process (cf. \cite{CsR}, page 30)
$$
|W(t_n)-W(n)|=O(n^{1/4+\varepsilon})\qquad {\rm a.s.}
$$
as $n\to\infty$, we have (\ref{inv}) of Theorem 1.1. $\Box$

\renewcommand{\thesection}{\arabic{section}.}
\section{Proof of Theorem 1.2}
\renewcommand{\thesection}{\arabic{section}} \setcounter{equation}{0}
\setcounter{theorem}{0} \setcounter{lemma}{0}

For $R>0$ integer define
\begin{eqnarray*}
\kappa_1&:=&\min\{t\geq 0:\, Y_\nu(t)=R\},\\
\delta_1&:=&\min\{t\geq \kappa_1:\, Y_\nu(t)\notin (R-1,R+1)\},\\
\kappa_i&:=&\min\{t\geq \delta_{i-1}:\, Y_\nu(t)=R\},\\
\delta_i&:=&\min\{t\geq \kappa_i:\, Y_\nu(t)\notin (R-1,R+1)\},\\
\kappa^*&:=&\max\{t\geq 0\, :Y_{\nu}(t)=R\},
\end{eqnarray*}
$i=2,3,\ldots$

Consider the local times at $R$ of the Bessel process during
excursions around $R$, i.e. let
$$
\zeta_i:=\eta(R,\delta_i)-\eta(R,\kappa_i),\quad i=1,2,\ldots,
$$
$$\tilde\zeta:=\eta(R,\infty)-\eta(R,\kappa^*).$$
We have
$$
\eta(R,\infty)=\sum_{i=1}^{\xi(R,\infty)-1}\zeta_i+\tilde\zeta.
$$
\begin{lemma}
\begg \be\left(e^{\lambda \eta(R,\infty)} \right)=\frac{p^*_R\,
\varphi(\lambda)}{1-q_R^*\,\varphi(\lambda)},
\endd
where \begg
p_R^*=\frac{A_R}{A_R+B_R}\,\frac{(R+1)^{2\nu}-R^{2\nu}}{(R+1)^{2\nu}},\quad
q_R^*=1-p_R^*,
\endd

\begg \varphi(\lambda)=\frac{\nu (A_R+B_R)}{\nu (A_R+B_R)-\lambda
R^{2\nu+1}A_R B_R}, \label{filam}
\endd
and \begg A_R=(R-1) ^{-2\nu}-R ^{-2\nu},\qquad
B_R=R^{-2\nu}-(R+1) ^{-2\nu}.
\endd
\end{lemma}

\noindent {\bf Proof:} By (\cite{BS}, p. 395, 3.3.2) $\zeta_i$
are i.i.d. random variables having exponential distribution with
moment generating function $\varphi(\lambda)$ given in
(\ref{filam}). Moreover, it is obvious that $\tilde\zeta$ is
independent from $\sum_{i=1}^{\xi(R,\infty)-1}\zeta_i.$
Furthermore, $\tilde\zeta$ is the  local time of $R$ under the
condition that starting from  $R, $ $Y_\nu(t)$ will reach $R+1$
 before $R-1.$  Hence its distribution can be calculated  from formula 3.3.5(b)
 of \cite{BS}, and its moment generating  function   happens to be equal to
$\varphi(\lambda)$  of (\ref{filam}). $\Box$

We can see
$$
\theta:=\be(\zeta_i)=\be(\tilde\zeta)=\frac{\nu(A_R+B_R)}
{R^{2\nu+1}A_R B_R}=1+O\left(\frac1R\right),\quad R\to\infty.
$$

$$\bp(|\eta(R,\infty)-\xi(R,\infty)|\geq u)=
\bp\left(\left|\sum_{i=1}^{\xi(R,\infty)-1}(\zeta_i-\theta)+
\tilde\zeta-\theta\right|\geq u\right)$$
$$\leq \bp(\xi(R.\infty)>N)+\bp \left(\max_{k\leq
N}\left|\sum_{i=1}^k(\zeta_i-\theta)\right|\geq u\right)
$$

$$
\leq (q^*_R)^N+ e^{-\lambda u}\left(\left(\frac{e^{\lambda
\theta}}{1+\lambda \theta}\right)^N+ \left(\frac{e^{-\lambda
\theta}}{1-\lambda \theta}\right)^N\right).
$$
In the above calculation we used the common moment generating
function (\ref{filam}) of $\zeta_i $ and $\tilde\zeta$, the exact
distribution of $\xi(R,\infty)$ (see (\ref{locelo})) and the exponential
Kolmogorov inequality. Estimating the above expression with
standard methods and selecting
$$N=CR\log R, \quad u=CR^{1/2} \log R,\quad \lambda=\frac{u}{\theta^2 N}$$
we conclude that
$$\bp(|\eta(R,\infty)-\xi(R,\infty)|\geq CR^{1/2} \log R )
\leq C_1 \exp{\left(-\frac{C\log R}{2 \theta}\right)}.$$ With a big
enough $C$ the right hand side of the above inequality is summable
in $R,$ hence Theorem 1.2 follows by the Borel-Cantelli lemma.
$\Box$

\renewcommand{\thesection}{\arabic{section}.}
\section{Proof of Theorem 1.3}
\renewcommand{\thesection}{\arabic{section}} \setcounter{equation}{0}
\setcounter{theorem}{0} \setcounter{lemma}{0}

Let $p_j^{(1)}$ and $p_j^{(2)}$ as in Theorem 1.3.
Define the two-dimensional Markov chain $(X_n^{(1)},X_n^{(2)})$ as
follows. If $p_j^{(1)}\geq p_k^{(2)}$, then let
\begin{eqnarray*}
\bp\left((X_{n+1}^{(1)},X_{n+1}^{(2)})=(j+1,k+1)\mid
(X_n^{(1)},X_n^{(2)})=(j,k)\right)&=&\frac12+p_k^{(2)}\\
\bp\left((X_{n+1}^{(1)},X_{n+1}^{(2)})=(j+1,k-1)\mid
(X_n^{(1)},X_n^{(2)})=(j,k)\right)&=&p_j^{(1)}-p_k^{(2)}\\
\bp\left((X_{n+1}^{(1)},X_{n+1}^{(2)})=(j-1,k-1)\mid
(X_n^{(1)},X_n^{(2)})=(j,k)\right)&=&\frac12-p_j^{(1)}.
\end{eqnarray*}
If, however $p_j^{(1)}\leq p_k^{(2)}$, then let
\begin{eqnarray*}
\bp\left((X_{n+1}^{(1)},X_{n+1}^{(2)})=(j+1,k+1)\mid
(X_n^{(1)},X_n^{(2)})=(j,k)\right)&=&\frac12+p_j^{(1)}\\
\bp\left((X_{n+1}^{(1)},X_{n+1}^{(2)})=(j-1,k+1)\mid
(X_n^{(1)},X_n^{(2)})=(j,k)\right)&=&p_k^{(2)}-p_j^{(1)}\\
\bp\left((X_{n+1}^{(1)},X_{n+1}^{(2)})=(j-1,k-1)\mid
(X_n^{(1)},X_n^{(2)})=(j,k)\right)&=&\frac12-p_k^{(2)}.
\end{eqnarray*}

Then it can be easily seen that $X_n^{(1)}$ and $X_n^{(2)}$ are
two NN random walks as desired. Consider the following 4 cases.
\begin{itemize}
\item{(i)} $p_j^{(1)}\leq p_k^{(2)}$, $j\leq k$,
\item{(ii)} $p_j^{(1)}\leq p_k^{(2)}$, $j\geq k$,
\item{(iii)} $p_j^{(1)}\geq p_k^{(2)}$, $j\leq k$,
\item{(iv)} $p_j^{(1)}\geq p_k^{(2)}$, $j\geq k$.
\end{itemize}

In case (i) from (\ref{pj1}) and (\ref{pj2}) we obtain
$$
\frac{B}{4j}-\frac{C}{j^\gamma}\leq
\frac{B}{4k}+\frac{C}{k^\gamma} \leq
\frac{B}{4k}+\frac{C}{kj^{\gamma-1}},
$$
implying
$$
k-j\leq \frac{2Cj^{2-\gamma}}{B/4-Cj^{1-\gamma}}=
O(j^{2-\gamma})
$$
if $j\to\infty$. So in this case if $X_n^{(1)}=j$ and
$X_n^{(2)}=k$, then we have
$$
0\leq X_n^{(2)}-X_n^{(1)}=O((X_n^{(1)})^{2-\gamma})
$$
if $n\to\infty$.

In case (ii) either
$X_{n+1}^{(1)}-X_{n+1}^{(2)}=X_{n}^{(1)}-X_{n}^{(2)}$, or
$X_{n+1}^{(1)}-X_{n+1}^{(2)}=X_{n}^{(1)}-X_{n}^{(2)}-2$, so that
we have
$$
-2\leq X_{n+1}^{(1)}-X_{n+1}^{(2)}\leq X_{n}^{(1)}-X_{n}^{(2)}.
$$

Similar procedure shows that in case (iii)
$$
-2\leq X_{n+1}^{(2)}-X_{n+1}^{(1)}\leq X_{n}^{(2)}-X_{n}^{(1)}
$$
and in case (iv)
$$
0\leq X_n^{(1)}-X_n^{(2)}=O((X_n^{(2)})^{2-\gamma}).
$$
Hence Theorem 1.3 follows from the law of the iterated logarithm for
$X_n^{(i)}$ (cf. \cite{BRS}). $\Box$

\renewcommand{\thesection}{\arabic{section}.}
\section{Strong theorems}
\renewcommand{\thesection}{\arabic{section}} \setcounter{equation}{0}
\setcounter{theorem}{0} \setcounter{lemma}{0}

As usual, applying Theorem 1.1 and Theorem 1.3, we can give limit
results valid for one of the processes to the other process
involved.

In this section we denote $Y_\nu(t)$ by $Y(t)$ and define the
following related processes.

$$
M(t):=\max_{0\leq s\leq t}Y(s), \qquad Q_n:=\max_{1\leq k\leq n}X_k.
$$

The future infimums are defined as
$$
I(t):=\inf_{s\geq t} Y(s), \qquad J_n:=\inf_{k\geq n}X_k.
$$

Escape processes are defined by
$$
A(t):=\sup\{s:\, Y(s)\leq t\}, \qquad G_n:=\sup\{k:\, X_k\leq n\}.
$$

Laws of the iterated logarithm are known for  Bessel processes
(cf. \cite{BS}) and NN random walks (cf. \cite{BRS}) as well.
Upper class results for Bessel process read as follows (cf. Orey and
Pruitt \cite{OP} for integral $d$, and Pardo \cite{PA} for the case of
positive self-similar Markov processes).

\medskip\noindent
{\bf Theorem F}: {\it Let $a(t)$ be a non-decreasing non-negative
continuous function. Then for $\nu\geq 0$
$$
\displaystyle{ t^{1/2}a(t)\in {\rm UUC}(Y(t))\qquad  if \,\, and
\,\, only \,\, if
\qquad\int_1^\infty\frac{(a(x))^{2\nu+2}}{x}e^{-a^2(x)/2}\,
dx<\infty.}
$$}

Now Theorems 1.1,  1.3 and Theorems E and F together imply the
following result.
\begin{theorem}
Let $\{X_n\}$ be an NN random walk with $p_R$ satisfying
$$
p_R=\frac{B}{4R}+O\left(\frac1{R^{1+\delta}}\right),\quad
R\to\infty
$$
with $B>1$ and for some $\delta>0$. Let furthermore $a(t)$ be
a  non-decreasing non-negative function. Then
$$
\displaystyle{n^{1/2}a(n)\in {\rm UUC}(X_n) \qquad if \,\, and
\,\, only \,\, if \qquad
\sum_{k=1}^\infty\frac{(a(k))^{B+1}}{k}e^{-a^2(k)/2}<\infty.}
$$

If $b(t)$ is a non-increasing non-negative function, then
$$
\displaystyle{n^{1/2}b(n)\in {\rm LLC}(X_n) \qquad if \,\, and
\,\, only \,\, if \qquad
 \sum_{k=1}^\infty (b(2^k))^{B-1}<\infty.}
$$
\end{theorem}

Next we prove the following invariance principles for the
processes defined above.

\begin{theorem}
Let $Y(t)$ and $X_n$ as in {\rm Theorem 1.1.} Then for any
$\varepsilon>0$ we have
\begg
|M(n)-Q_n|=O(n^{1/4+\varepsilon})\quad {\rm a.s.}
\label{invmq}
\endd
and
\begg
|I(n)-J_n|=O(n^{1/4+\varepsilon})\quad {\rm a.s.}
\label{invij}
\endd
\end{theorem}

\noindent {\bf Proof:} Define $\tilde s, s^*, \tilde k, k^*$ by
$$
Y(\tilde s)=M(n),\quad Y(s^*)=I(n), \quad X_{\tilde k}=Q_n, \quad
X_{k^*}=J_n.
$$
Then as $n\to\infty$, we have almost surely
$$
Q_n-M(n)=X_{\tilde k}-Y(\tilde s)\leq X_{\tilde k}-Y(\tilde k)
=O(n^{1/4+\varepsilon})
$$
and
$$
M(n)-Q_n=Y(\tilde s)-X_{\tilde k}=Y(\tilde s)-Y([\tilde s])
-(X_{[\tilde s]}-Y([\tilde s]))+X_{[\tilde s]}-X_{\tilde k}
$$
$$
\leq Y(\tilde s)-Y([\tilde s])
-(X_{[\tilde s]}-Y([\tilde s])=Y(\tilde s)-Y([\tilde s])
+O(n^{1/4+\varepsilon})
$$
By (\ref{dife}) and recalling the results on the increments of the
Wiener process (see \cite{CsR} page 30) we get
$$
Y(\tilde s)-Y([\tilde s])=W(\tilde s)-W([\tilde s])+
\int_{[\tilde s]}^{\tilde s}\frac{2\nu+1}{2Y(s)}\, ds
$$
$$
\leq \sup_{0\leq t\leq n}\sup_{0\leq s\leq 1}|W(t+s)-W(t)|
+\frac{2\nu+1}{2}\max_{[\tilde s]\leq t\leq \tilde s}\frac{1}{Y(t)}
=O(\log n) \quad{\rm a.s.},
$$
since $Y(t)$ in the interval $([\tilde s],\tilde s)$ is bounded away
from zero. Hence (\ref{invmq}) follows.

To show (\ref{invij}), note that $n\leq s^*\leq n^{1+\alpha}$ and $n\leq
k^*\leq n^{1+\alpha}$ for any $\alpha>0$ almost surely for all large
$n$.
Then as $n\to\infty$
$$
I(n)-J_n\leq Y(k^*)-X_{k^*}=O((k^*)^{1/4+\varepsilon})=
O(n^{(1+\alpha)(1/4+\varepsilon)})\quad {\rm a.s.}
$$
On the other hand,
$$
J_n-I(n)\leq X_{k^*}-Y([s^*])+Y([s^*])-Y(s^*)=
O(n^{(1+\alpha)(1/4+\varepsilon)})+Y([s^*])-Y(s^*).
$$
By (\ref{dife}), taking into account that when applying this
formula the integral contribution is negative, and recalling again the
results on the increments of the Wiener process, we get
$$
Y([s^*])-Y(s^*)\leq W([s^*])-W(s^*)\leq \sup_{0\leq t\leq
n^{1+\alpha}}\sup_{0\leq s\leq 1}|W(t+s)-W(t)|=O(\log n)\quad
{\rm a.s.}
$$
as $n\to\infty$. Hence
$$
|I(n)-J_n|=O(n^{(1+\alpha)(1/4+\varepsilon)})\quad {\rm a.s.}
$$
Since $\alpha>0$ and $\varepsilon>0$ are arbitrary, (\ref{invij})
follows. This completes the proof of Theorem 6.2.
$\Box$

\begin{theorem}
Let $X_n^{(1)}$ and $X_n^{(2)}$ as in {\rm Theorem 1.3} and let
$Q_n^{(1)}$ and $Q_n^{(2)}$ be the corresponding maximums, while let
$J_n^{(1)}$ and $J_n^{(2)}$ be the corresponding future infimum
processes. Then for any $\varepsilon>0$, as $n\to\infty$ we have
\begg
|Q_n^{(1)}-Q_n^{(2)}|=O(n^{1-\gamma/2+\varepsilon})\quad {\rm a.s.}
\label{invq1q2}
\endd
and
\begg
|J_n^{(1)}-J_n^{(2)}|=O(n^{1-\gamma/2+\varepsilon})\quad {\rm a.s.}
\label{invj1j2}
\endd
\end{theorem}

\noindent {\bf Proof:} Define $\tilde k_i, k_i^*,\, i=1,2$ by
$$
X^{(i)}_{\tilde k_i}=Q_n^{(i)}, \quad
X_{k_i^*}^{(i)}=J_n^{(i)}.
$$
Then
$$
|Q_n^{(1)}-Q_n^{(2)}|\leq \max(X_{\tilde k_1}^{(1)}-X_{\tilde
k_1}^{(2)}, X_{\tilde k_2}^{(1)}-X_{\tilde k_2}^{(2)})
=O((n\log\log n)^{1-\gamma/2})\quad {\rm a.s.},
$$
proving (\ref{invq1q2}).

Moreover, for any $\alpha>0$, $n\leq k_i^*\leq n^{1+\alpha}$ almost
surely for large $n$, hence we have
$$
|J_n^{(1)}-J_n^{(2)}|\leq \max(X_{k_1^*}^{(1)}-X_{k_1^*}^{(2)},
X_{k_2^*}^{(1)}-X_{k_2^*}^{(2)})
=O((n\log\log n)^{(1+\alpha)(1-\gamma/2)})\quad {\rm a.s.}
$$
Since $\alpha$ is arbitrary, (\ref{invj1j2}) follows.

This completes the proof of Theorem 6.3. $\Box$

\medskip

Khoshnevisan et al. \cite{KLS} (for $I(t)$ and $A(t)$), Adelman
and Shi \cite{AS}, and Shi \cite{Shi} (for $Y(t)-I(t)$) proved
the following upper and lower class results.

\medskip\noindent
{\bf Theorem G}: {\it Let $\varphi(t)$ be a non-increasing, and
$\psi(t)$ be a non-decreasing function, both non-negative. Then for}
$\nu>0$
\begin{itemize}
 \item $\displaystyle{t^{1/2}\psi(t)\in {\rm UUC}(I(t))}$
\qquad {\it if and only if} \qquad
$\displaystyle{\int_1^\infty\frac{(\psi(x))^{2\nu}}{x}e^{-\psi^2(x)/2}\,
dx<\infty,}$
\item
$\displaystyle{t^2\varphi(t)\in {\rm LLC}(A(t))} $ \qquad {\it if
and only if} \qquad $
\displaystyle{\int_1^\infty\frac1{x\varphi^{\nu}(x)}e^{-1/2\varphi(x)}\,
dx<\infty. }$
\item
$ \displaystyle{t^{1/2}\psi(t)\in {\rm UUC}(Y(t)-I(t))} $ \qquad
{\it if and only if} \qquad $
\displaystyle{\int_1^\infty\frac1{x\psi^{2\nu-2}(x)}e^{-\psi^2(x)/2}\,
dx<\infty,} $
\end{itemize}
\medskip\noindent
{\bf Theorem H}: { \it Let $\rho(t)>0$ be such that
$(\log\rho(t))/\log t$ is non-decreasing.} Then
\begin{itemize}
\item $\displaystyle{1/\rho(t)\in {\rm LLC}(M(t)-I(t)) }$ \qquad
{\it if and only if} \qquad $ \displaystyle{\int_1^\infty
\frac{dx}{x\log\rho(x)}<\infty.}$
\end{itemize}

Taking into account that $J_n$ and $G_n$ are inverses of each other,
immediate consequences of Theorems F,\, G,\, H, Theorems 6.2 and 6.3 are
the following upper and lower class results.

\begin{theorem}
Let $X_n$ be as in {\rm Theorem 6.1} and let $\varphi(t)$ be a
non-increasing and $\psi(t)$ be a non-decresing function, both
non-negative. Then
\begin{itemize}
\item $ \displaystyle{n^{1/2}\psi(n)\in {\rm UUC}(J_n)}$
 \qquad if and only if \qquad
$\displaystyle{\sum_{k=1}^\infty\frac{(\psi(k))^{B-1}}{k}
e^{-\psi^2(k)/2}<\infty,}$
\item
$\displaystyle{n^2\varphi(n)\in {\rm LLC}(G_n) }$ \qquad if and only
if \qquad
$\displaystyle{\sum_{k=1}^\infty\frac1{k\varphi^{(B-1)/2}(k)}e^{-1/2\varphi(k)}
<\infty.} $
\item
$\displaystyle{ n^{1/2}\psi(n)\in {\rm UUC}(X_n-J_n)} $ \qquad if
and only if \qquad $
\displaystyle{\sum_{k=1}^\infty\frac1{k\psi^{B-3}(k)}e^{-\psi^2(k)/2}<\infty,}
$
\end{itemize}
\end{theorem}
\begin{theorem} Let $\rho(t)>0$ be such that
$(\log\rho(t))/\log t$ is non-decreasing.
\begin{itemize}
\item
$\displaystyle{ 1/\rho(n)\in {\rm LLC}(Q_n-J_n) }$ \qquad if and
only if \qquad $ \displaystyle{\sum_{k=2}^\infty \frac1{k\log
\rho(k)}<\infty.} $
\end{itemize}
\end{theorem}

\renewcommand{\thesection}{\arabic{section}.}
\section{Local time}
\renewcommand{\thesection}{\arabic{section}} \setcounter{equation}{0}
\setcounter{theorem}{0} \setcounter{lemma}{0}

We will need the following result from Yor \cite{Yor}, page 52.

\noindent
{\bf Theorem J:} {\it For the local time of a Bessel process of order
$\nu$ we have}
$$\eta(R,\infty)\eqe
(2\nu)^{-1}R^{1-2\nu}Y_0^2(R^{\,2\nu}),
$$
{\it where $Y_0$ is a two-dimensional Bessel process and $\eqe$ means
equality in distribution.}

Hence applying Theorem F for $\nu=0$, we get

\medskip\noindent
{\bf Theorem K}: {\it If $f(x)$ is non-decreasing, non-negative
function, then}
\begin{itemize}
\item $\displaystyle{ Rf(R)\in {\rm UUC}(\eta(R,\infty))}$
\qquad {\it if and only if} \qquad $\displaystyle{\int_1^\infty
\frac{f(x)}{x}e^{-\nu f(x)}\, dx<\infty}. $
\end{itemize}

From this and Theorem 1.2 we get the following result.
\begin{theorem}{\it If $f(x)$ is non-decreasing, non-negative function,
then}
\begin{itemize}
\item $\displaystyle{
Rf(R)\in {\rm UUC}(\xi(R,\infty))}$ \qquad if and only if\qquad
$\displaystyle{ \sum_{k=1}^\infty \frac{f(k)}{ k} e^{-\nu
f(k)}<\infty.} $
\end{itemize}
\end{theorem}

In \cite{CSFR} we proved the following result.

\medskip\noindent
{\bf Theorem L:} {\it Let $\displaystyle{
p_R=\frac{B}{4R}+O\left(\frac{1}{R^\gamma}\right)}$ with $ B>1$, and
$\gamma>1$. Then with probability $1$ there exist infinitely many
$R$ for which
$$\xi(R+j,\infty)=1$$
for each $j=0,1,2,\ldots,[\log\log R/ \log 2]$. Moreover, with
probability $1$ for each $R$ large enough and $\varepsilon>0$
there exists an
$$R\leq S\leq {\frac{(1+\varepsilon)\log\log R}{\log 2}}$$
such that
$$\xi(S,\infty)>1.$$
}

\noindent {\bf Remark 1:} In fact in \cite{CSFR} we proved this
result in the case when $p_R= B/{4R}$  but the same proof works also
in the case of Theorem L.

This theorem applies e.g. for the case when $p_R$ is given by
(\ref{pr}), which in turn, gives the following result for the Bessel
process.

Let
\begin{itemize}
\item[(i)] $\kappa(R):=\inf\{t:\ Y_\nu(t)=R\}$,
\item[(ii)]
$\kappa^*(R):=\sup\{t:\ Y_\nu(t)=R\}$,
\item[(iii)] $\Psi(R)$ be the largest
integer for which the event
$$A(R)=\bigcap_{j=-1}^{\Psi(R)}\{\kappa^*(R+j)<\kappa(R+j+1)\}$$
occurs.
\end{itemize}
$A(R)$ means that $Y_\nu(t)$ moves from $R$ to $R+1$ before
returning to $R-1$, it goes from $R+1$ to $R+2$ before returning
to $R$, $\ldots$ and also from $R+\Psi(R)$ to $R+\Psi(R)+1$ and it
never returns to $R+\Psi(R)-1$. We say that the process
$Y_\nu(t)$ escapes  through $(R,R+\Psi(R))$ with large velocity.

\begin{theorem}
$$\limsup_{R\to\infty}{\frac{\Psi(R)}{\log\log R}}=
{\frac{1}{\log 2}}\quad {\rm a.s.}$$
\end{theorem}
\vspace{2ex}\noindent {\bf Remark 2:} The statement of Theorem 7.2
(for integral $d=2\nu+2$) was formulated in \cite{R05}, p. 291 as a
Conjecture.

\end{document}